\documentclass[12pt]{article}
\usepackage[utf8]{inputenc}
\usepackage{mathrsfs}
\usepackage{fullpage}
\usepackage{amssymb}
\usepackage{amsmath}
\usepackage{epsfig}
\usepackage{algorithm}
\usepackage{algorithmic}
\usepackage{latexsym}
\usepackage{amsmath}
\usepackage{amsfonts}
\usepackage{graphicx}
\usepackage{authblk}
\usepackage[normalem]{ulem}
\usepackage{array}
\usepackage{colortbl}
\usepackage{bm}
\usepackage{color}\usepackage[dvipsnames]{xcolor}
\usepackage{tikz}\usetikzlibrary{patterns}
\usepackage{chessboard}
\usepackage{skak}
\usepackage{array,multirow}
\usepackage[
  bookmarks=false, 
  colorlinks,
  citecolor=brown!70!black,
  linkcolor=brown!80!black,
  urlcolor=blue!70!black,
]{hyperref}
\usepackage{relsize,exscale}

\newtheorem{defn}{Definition}
\newtheorem{cor}{Corollary}
\newtheorem{thm}{Theorem}

\newcommand{\comm}[1]{}
 	\definecolor{lightlightgray}{rgb}{0.93, 0.93, 0.93}
 		\definecolor{llightgray}{rgb}{0.87, 0.87, 0.87}

\newcolumntype{C}[1]{>{\centering\arraybackslash }b{#1}}

\def\A{0.5cm}

 \def\C{2.5cm}
 
 \def\E{4.5cm}

 \def\G{6.5cm}
 
 \def\I{8.5cm}

 \def\K{10.5cm}
 \def\LL{11.5cm}
\def\M{12.5cm}
 \def\N{13.5cm}
 \def\O{14.5cm}\def\Q{16.5cm}\def\R{17.5cm}\def\S{18.5cm}
\def\U{20.5cm}\def\V{21.5cm}\def\W{22.5cm}\def\Y{24.5cm}\def\ZZ{26.5cm}\def\Za{28.5cm}\def\Zb{30.5cm}\def\Zc{32.5cm}\def\Zd{34.5cm}\def\Ze{36.5cm}

\def\sizePoint{2.6pt}

\newcommand{\point}[2]{\fill (canvas cs:x=#1,y=#2) circle (\sizePoint); }

\newcommand\oeis[1]{\href{https://oeis.org/#1}{#1}}

\title{Knight's paths towards Catalan numbers}

\author[1]{Jean-Luc Baril}
\author[2]{Jos{\'e} L. Ram\'irez}
\affil[1]{\rm LIB, Universit\'e de Bourgogne Franche-Comt\'e \protect\\
  B.P. 47 870, 21078 Dijon Cedex France\protect\\
   {\tt E-mail: barjl@u-bourgogne.fr
   }
}

\affil[2]{\rm Departamento de Matem{\'a}ticas, Universidad  Nacional de Colombia\protect\\
 Bogot{\'a}, Colombia\protect\\
   {\tt E-mail: jlramirezr@unal.edu.co
   }
}
\date{\today}
\providecommand{\keywords}[1]{\textit{Keywords:} #1}
\providecommand{\subjclass}[1]{\textit{MSC:} #1}

\begin{document}

\maketitle 

\begin{abstract} We provide enumerating results for partial knight's paths of a given size. We prove algebraically that zigzag knight's paths of a given size ending on the $x$-axis are enumerated by the generalized Catalan numbers, and  we give a constructive bijection with peakless Motzkin paths of a given length. After enumerating partial knight's paths of a given length, we prove that zigzag knight's paths of a given length ending on the $x$-axis are counted by the Catalan numbers. Finally, we give a constructive bijection with  Dyck paths of a given length.
\end{abstract}


\section{Introduction}
 For several decades  lattice paths have been  widely studied in the literature. Their links with  many problems of various domains as computer science, biology and physics~\cite{Sta}, give
them an important place in combinatorics. For instance, they have very tight connections with RNA structures, pattern avoiding permutations, directed animals, and so on~\cite{Bar1,Dos,Knu,Kra,Sta}.  In combinatorics, a classical problem consists in  enumerating some special kinds of paths with respect to the length and other statistics (see for instance \cite{Ban,Barc,Bar,Deu,Flo, Man1,Mer,Sap,Sun}).  Certainly,  Dyck paths are the most known. These are lattice paths in $\Bbb N^2=\{0,1,2,\ldots\}^2$ starting at the origin $(0,0)$, ending on the $x$-axis, consisting of  steps lying in $\{(1,1), (1,-1)\}$. They  are  enumerated with respect to the number of steps  by the  famous  Catalan  numbers (see \href{https://oeis.org/A000108}{A000108} in Sloane's On-line Encyclopedia  of  Integer Sequences~\cite{Sloa}). 

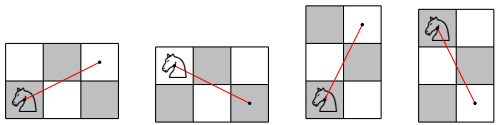
\begin{figure}[H]
 \begin{center}
        \begin{tikzpicture}[scale=0.25]
           
            \fill[lightgray] (0,0) -- (2,0) -- (2,2) -- (0,2) -- cycle;
                  \fill[lightgray] (2,2) -- (2,4) -- (4,4) -- (4,2) -- cycle;
                       \fill[lightgray] (4,0) -- (4,2) -- (6,2) -- (6,0) -- cycle;
                       \draw (0,0)-- (6,0);
                        \node at (1,1) {\knight} ;
             \draw (0,2)-- (6,2);
              \draw (0,4)-- (6,4);
              \draw (0,0)-- (0,4);
               \draw (2,0)-- (2,4);
                \draw (4,0)-- (4,4);
                 \draw (6,0)-- (6,4);
                
                 \draw[red] (1,1) -- (5,3);
                 \point{1cm}{1cm}
                 \point{5cm}{3cm}
     \end{tikzpicture}\quad
  \begin{tikzpicture}[scale=0.25]
   \fill[lightgray] (0,0) -- (2,0) -- (2,2) -- (0,2) -- cycle;
                  \fill[lightgray] (2,2) -- (2,4) -- (4,4) -- (4,2) -- cycle;
                       \fill[lightgray] (4,0) -- (4,2) -- (6,2) -- (6,0) -- cycle;
                         \node at (1,3) {\knight} ;
            \draw (0,0)-- (6,0);
             \draw (0,2)-- (6,2);
              \draw (0,4)-- (6,4);
              \draw (0,0)-- (0,4);
               \draw (2,0)-- (2,4);
                \draw (4,0)-- (4,4);
                 \draw (6,0)-- (6,4);
                 \draw[red] (1,3) -- (5,1);
                 \point{1cm}{3cm}
                 \point{5cm}{1cm}
     \end{tikzpicture}\quad
\begin{tikzpicture}[scale=0.25]
 \fill[lightgray] (0,0) -- (2,0) -- (2,2) -- (0,2) -- cycle;
                  \fill[lightgray] (2,2) -- (2,4) -- (4,4) -- (4,2) -- cycle;
                       \fill[lightgray] (0,4) -- (2,4) -- (2,6) -- (0,6) -- cycle;
                         \node at (1,1) {\knight} ;
            \draw (0,0)-- (4,0);
             \draw (0,2)-- (4,2);
              \draw (0,4)-- (4,4);
              \draw (0,6)-- (4,6);
               \draw (0,0)-- (0,6);
                \draw (2,0)-- (2,6);
                 \draw (4,0)-- (4,6);
                 \draw[red] (1,1) -- (3,5);
                 \point{1cm}{1cm}
                 \point{3cm}{5cm}
     \end{tikzpicture}\quad
     \begin{tikzpicture}[scale=0.25]
      \fill[lightgray] (0,0) -- (2,0) -- (2,2) -- (0,2) -- cycle;
                  \fill[lightgray] (2,2) -- (2,4) -- (4,4) -- (4,2) -- cycle;
                       \fill[lightgray] (0,4) -- (2,4) -- (2,6) -- (0,6) -- cycle;
                         \node at (1,5) {\knight} ;
            \draw (0,0)-- (4,0);
             \draw (0,2)-- (4,2);
              \draw (0,4)-- (4,4);
              \draw (0,6)-- (4,6);
               \draw (0,0)-- (0,6);
                \draw (2,0)-- (2,6);
                 \draw (4,0)-- (4,6);
                 \draw[red] (3,1) -- (1,5);
                 \point{3cm}{1cm}
                 \point{1cm}{5cm}
     \end{tikzpicture}
               \end{center}
         \caption{ The four possible right-moves  of a knight on a chessboard.}
         \label{fig1}
\end{figure}

In this work we deal with several classes of knight's paths (see \cite{Lab}), i.e.,  lattice paths in $\Bbb N^2$ starting at $(0,0)$, and consisting of right-moves  of a knight in the  game of chess (a  right-move of a knight is a jump from left to right from one corner of any two-by-three rectangle to the opposite corner on a chessboard). Figure \ref{fig1} shows the four possible right-moves of a knight. It is worth noticing  that knight's paths always start at the origin and  go to the right at each step. Therefore, this work is not related to the well known mathematical problem about the walk of a knight on a chessboard, namely the knight's tour problem \cite{Eul}, which consists in finding a sequence of moves on a chessboard such that the knight visits every square exactly once.

\begin{figure}[h]
 \begin{center}
        \begin{tikzpicture}[scale=0.15]
            \draw (\A,\A)-- (38,\A);
             \draw[dashed,line width=0.1mm] (\A,\E)-- (\Ze,\E);
              \draw[dashed,line width=0.1mm] (\A,\C)-- (\Ze,\C);
               \draw[dashed,line width=0.1mm] (\A,\G)-- (\Ze,\G);
               \draw[dashed,line width=0.1mm] (\A,\I)-- (\Ze,\I);
              \draw[dashed,line width=0.1mm] (\A,\K)-- (\Ze,\K);
               \draw[dashed,line width=0.1mm] (\A,\M)-- (\Ze,\M);
            \draw (\A,\A) -- (\A,\O);
             \draw[dashed,line width=0.1mm] (\C,\A) -- (\C,\M);\draw[dashed,line width=0.1mm] (\E,\A) -- (\E,\M);\draw[dashed,line width=0.1mm] (\G,\A) -- (\G,\M);
             \draw[dashed,line width=0.1mm] (\I,\A) -- (\I,\M);\draw[dashed,line width=0.1mm] (\K,\A) -- (\K,\M);\draw[dashed,line width=0.1mm] (\M,\A) -- (\M,\M);
             \draw[dashed,line width=0.1mm] (\O,\A) -- (\O,\M);\draw[dashed,line width=0.1mm] (\Q,\A) -- (\Q,\M);\draw[dashed,line width=0.1mm] (\S,\A) -- (\S,\M);
             \draw[dashed,line width=0.1mm] (\U,\A) -- (\U,\M);\draw[dashed,line width=0.1mm] (\W,\A) -- (\W,\M);\draw[dashed,line width=0.1mm] (\Y,\A) -- (\Y,\M);
             \draw[dashed,line width=0.1mm] (\ZZ,\A) -- (\ZZ,\M);
             \draw[dashed,line width=0.1mm] (\Za,\A) -- (\Za,\M);
             \draw[dashed,line width=0.1mm] (\Zb,\A) -- (\Zb,\M);
             \draw[dashed,line width=0.1mm] (\Zc,\A) -- (\Zc,\M);
             \draw[dashed,line width=0.1mm] (\Zd,\A) -- (\Zd,\M);
             \draw[dashed,line width=0.1mm] (\Ze,\A) -- (\Ze,\M);
            \draw[solid,line width=0.4mm] (\A,\A)--(\C,\E)  -- (\G,\C) -- (\I,\G) --(\K,\C)-- (\O,\E) -- (\S,\C) -- (\U,\G) -- (\Y,\E)  -- (\Za,\G) -- (\Zb,\C) -- (\Zd, \E)--(\Ze,\A);
         \end{tikzpicture}\qquad 
 \begin{tikzpicture}[scale=0.15]
            \draw (\A,\A)-- (38,\A);
             \draw[dashed,line width=0.1mm] (\A,\E)-- (\Ze,\E);
              \draw[dashed,line width=0.1mm] (\A,\C)-- (\Ze,\C);
               \draw[dashed,line width=0.1mm] (\A,\G)-- (\Ze,\G);
               \draw[dashed,line width=0.1mm] (\A,\I)-- (\Ze,\I);
              \draw[dashed,line width=0.1mm] (\A,\K)-- (\Ze,\K);
               \draw[dashed,line width=0.1mm] (\A,\M)-- (\Ze,\M);
            \draw (\A,\A) -- (\A,\O);
             \draw[dashed,line width=0.1mm] (\C,\A) -- (\C,\M);\draw[dashed,line width=0.1mm] (\E,\A) -- (\E,\M);\draw[dashed,line width=0.1mm] (\G,\A) -- (\G,\M);
             \draw[dashed,line width=0.1mm] (\I,\A) -- (\I,\M);\draw[dashed,line width=0.1mm] (\K,\A) -- (\K,\M);\draw[dashed,line width=0.1mm] (\M,\A) -- (\M,\M);
             \draw[dashed,line width=0.1mm] (\O,\A) -- (\O,\M);\draw[dashed,line width=0.1mm] (\Q,\A) -- (\Q,\M);\draw[dashed,line width=0.1mm] (\S,\A) -- (\S,\M);
             \draw[dashed,line width=0.1mm] (\U,\A) -- (\U,\M);\draw[dashed,line width=0.1mm] (\W,\A) -- (\W,\M);\draw[dashed,line width=0.1mm] (\Y,\A) -- (\Y,\M);
             \draw[dashed,line width=0.1mm] (\ZZ,\A) -- (\ZZ,\M);
             \draw[dashed,line width=0.1mm] (\Za,\A) -- (\Za,\M);
             \draw[dashed,line width=0.1mm] (\Zb,\A) -- (\Zb,\M);
             \draw[dashed,line width=0.1mm] (\Zc,\A) -- (\Zc,\M);
             \draw[dashed,line width=0.1mm] (\Zd,\A) -- (\Zd,\M);
             \draw[dashed,line width=0.1mm] (\Ze,\A) -- (\Ze,\M);
            \draw[solid,line width=0.4mm] (\A,\A)--(\C,\E)  -- (\G,\C) -- (\I,\G) --(\K,\C)-- (\O,\E) -- (\S,\C) -- (\U,\G) -- (\Y,\E)  -- (\Za,\G) -- (\Zb,\C) -- (\Zc, \G)--(\Ze,\E);
         \end{tikzpicture}
               \end{center}
         \caption{ The left and the right drawing show both a partial knight's path of size $18$ and length $12$. The left knight's path ends at height $0$, whereas the right one ends at height~$2$.}
         \label{fig2}
\end{figure}
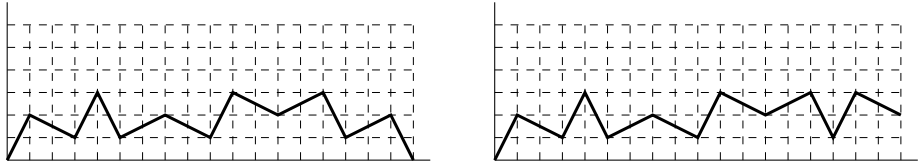

More formally, we define  knight's paths as follows.
\begin{defn}A {\it knight's path} is a lattice path in $\Bbb{N}^2$ starting at the origin, ending on the $x$-axis, and consisting of steps lying in $S=\{(1,2), (1,-2), (2,1),(2,-1)\}$. A {\it partial knight's path} is a prefix of a knight's path.
\end{defn}
The {\it length} of a path is the number of its steps, and the {\it size} of a path is the abscissa of its last point.  We denote by $\epsilon$ the empty path, i.e., the path of length $0$ (or equivalently of size $0$). The {\it height} of a point $(x,y)$ of a path is the ordinate $y$. We refer to Figure~\ref{fig2} for two examples of partial knight's paths where the heights of the last points are $0$ and $2$, respectively. Defining $E=(2,1)$, $\bar{E}=(2,-1)$, $N=(1,2)$ and $\bar{N}=(1,-2)$, these two paths can be written as $N\bar{E}N\bar{N}E\bar{E}N\bar{E}E\bar{N}E\bar{N}$ and $N\bar{E}N\bar{N}E\bar{E}N\bar{E}E\bar{N}N\bar{E}$. 
Knight's paths of a given size (ending on the $x$-axis) have been enumerated in $\cite{Lab}$. On the other hand, there are enumerative results for knight's paths of a given length which are in one-to-one correspondence with basketball walks in the first quarter plane (see for instance  \cite{Ayy,Bank, Bett, Bousq}). We will recall some of these results at the beginning of Sections 2 and 3.

In this paper, we will focus on  a special kind of knight's paths  that are enumerated by Catalan and generalized Catalan numbers (see \oeis{A004148} in \cite{Sloa} for a definition of the generalized Catalan numbers). More precisely, we will consider zigzag knight's paths defined as follows.
\begin{defn} A {\it zigzag knight's path} (or {\it zigzag path} for short) is a knight's path with the additional property that two consecutive steps cannot be in the same direction, i.e., two consecutive steps $(a,b)$ and $(c,d)$ lying in $S$ must satisfy $b\cdot d<0$. 
\end{defn}
In a zigzag path described as word in $N$, $\bar{N}$, $E$, and $\bar{E}$, barred and unbarred letters alternate. Notice that the two paths on Figure \ref{fig2} are zigzag paths.

\bigskip

{\bf Outline of the paper.} In Section 2, we provide enumerating results for partial knight's paths (prefixes of knight's paths) of a given size. We make a similar study for partial zigzag paths of a given size. We prove that zigzag paths of a given size ending on the $x$-axis are enumerated by the generalized Catalan numbers, and the total number of partial zigzag paths of a given size are enumerated as the symmetric secondary structures of RNA molecules with a given number of nucleotides. Next, we give a constructive one-to-one correspondence between zigzag paths of size $n$ and peakless Motzkin paths of length $n+1$. 

In Section 3,  we provide enumerating results for partial zigzag knight's paths of a given length.  We prove that zigzag paths of a given length ending on the $x$-axis are enumerated by the Catalan numbers, and the total number of partial zigzag paths of a given length are enumerated by the central binomial coefficient. Finally, we give a constructive one-to-one correspondence between zigzag paths of length $n$ and Dyck paths of length $2(n+1)$.


\section{Knight's paths of a given size}

In \cite{Lab}, Labelle and Yeh provide the generating function for the number of  knight's paths of size $n$ (ending on the $x$-axis):
$$A(z)=\frac{1+2z+\sqrt{1-4z+4z^2-4z^4}-\sqrt{2}\sqrt{1-4z^2-2z^4+(2z+1)\sqrt{1-4z+4z^2-4z^4}}}{4z^2}.$$
The first terms of the series expansion are
$$1+z^2+3z^4+2z^5+12z^6+14z^7+54z^8+86z^9+274z^{10}+528z^{11}+\ldots,$$ and the sequence of coefficients corresponds to \oeis{A005220} in \cite{Sloa}.  As mentioned in \cite{Lab} (see Theorem 2.13), the generating function $A(z)$ satisfies the functional equation 
$$z^4A(z)^4 - (2z^3 + z^2)A(z)^3 + (z^4+ 2z^2 + 2z)A(z)^2 - (2z+ 1)A(z) + 1=0.$$

\subsection{Partial knight's paths}
In this part, we count partial knight's paths of a given size, i.e., prefixes of knight's paths ending at a given abscissa. For $k\geq 0$, we consider the generating function $f_k=f_k(z)$ (resp. $g_k=g_k(z)$), where the coefficient of $z^n$ in the series expansion is the number of  partial knight's paths of size $n$  ending at height $k$ with an up-step $N$ or $E$, (resp. with a down-step  $\bar{N}$ or $\bar{E}$).  So, we easily obtain the following equations: 
\begin{equation}\begin{array}{l}
f_0=1, \mbox{ and } f_1=z^2(f_0+g_0),\\
f_k=z^2(f_{k-1}+g_{k-1})+z(f_{k-2}+g_{k-2}), \quad k\geq 2,\\
g_k=z^2(f_{k+1}+g_{k+1})+z(f_{k+2}+g_{k+2}), \quad k\geq 0.\\
\end{array}\label{eq1}
\end{equation}

Now, we introduce the bivariate generating functions
$$F(u,z)=\sum\limits_{k\geq 0} u^kf_k(z) \quad \mbox{and} \quad G(u,z)=\sum\limits_{k\geq 0} u^kg_k(z).$$
For short, we use the notation $F(u)$ and $G(u)$ for these functions, respectively.
Summing the recursions in (\ref{eq1}), we have:
\begin{align*}
F(u)&=1+z^2\sum\limits_{k\geq 1}u^k(f_{k-1}+g_{k-1}) + z\sum\limits_{k\geq 2}u^k(f_{k-2}+g_{k-2})\\
&=1+z^2u(F(u)+G(u))+zu^2(F(u)+G(u))\\
&=1+zu(z+u)(F(u)+G(u)),\\
G(u)&=z^2\sum\limits_{k\geq 0}u^k (f_{k+1}+ g_{k+1})+z\sum\limits_{k\geq 0}u^k (f_{k+2}+ g_{k+2})\\
&=\frac{z^2}{u}(F(u)-1+G(u)-g_0) + \frac{z}{u^2}(F(u)-1-uz^2(1+g_0)+G(u)-g_0-ug_1),
\end{align*}
where  $g_0=G(0)$ is the generating function for the number of non-empty knight's paths ending on the $x$-axis with respect to the size, i.e., $g_0=A(z)-1$, and $g_1$ is the generating function for the number of partial knight's paths ending at height $1$ with a down-step $\bar{N}$ or $\bar{E}$. Then, $g_1$ is the difference between the generating function $A_1(z)$ for the number of partial knight's paths ending at height 1, and the generating function for the number of partial knight's paths ending at height 1 with the step $E=(2,1)$, i.e. $g_1=A_1(z)-z^2(g_0+1)$ where $A_1(z)$ is given in \cite{Lab} (see Lemma 2.11) as follows:
$$A_1(z)=z^2A(z)^2(1-zA(z))^{-1}.$$

Solving the above functional equations, we deduce
$$F(u)={\frac {uz^2A(z)(u^2z+uz^2+u+z)+u^2z^2A_1(z)(u+z)+uz^2-u^2+z}{{u}^{4}z+{u}^{3}{z}^{2
}+u{z}^{2}-{u}^{2}+z}},$$

$$G(u)=-{\frac { \left(
A(z)(u^3z^2+u^2z^3+u^2z+uz^2-uz-1)+A_1(z)(u^3z+u^2z^2-u)+zu+1
\right) z}{{u}^{4}z+{u}^{3}{z}^
{2}+u{z}^{2}-{u}^{2}+z}},$$

 $$F(u)+G(u)={\frac {z(uz+1)A(z)+uzA_1(z)-{u}^{2}}{{u}^{4}z+{u}^{3}{z}^{2}+u{z}^{2}-{u}
^{2}+z}}
.$$
The first terms of the series expansion of $F(u)+G(u)$ are
\begin{multline*}
1 + u^2 z + (1 + u + u^4) z^2 + (u + 2 u^2 + 2 u^3 + u^6) z^3 + (3 + 3u + u^2 + 2 u^3 + 3 u^4 + 3 u^5 + u^8) z^4 \\+ (2 + 4 u  + 9 u^2 + 
    8 u^3 + 3 u^4 + 3 u^5 + 4 u^6 + 4 u^7 + u^{10}) z^5+\ldots
\end{multline*}

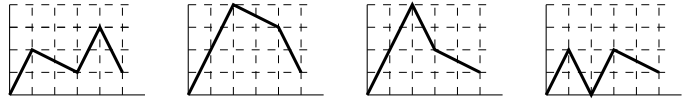
\begin{figure}[h]
 \begin{center}
        \begin{tikzpicture}[scale=0.15]
            \draw (\A,\A)-- (\M,\A);
             \draw[dashed,line width=0.1mm] (\A,\E)-- (\LL,\E);
              \draw[dashed,line width=0.1mm] (\A,\C)-- (\LL,\C);
               \draw[dashed,line width=0.1mm] (\A,\G)-- (\LL,\G);
            \draw[dashed,line width=0.1mm] (\A,\G)-- (\LL,\G);
             \draw[dashed,line width=0.1mm] (\A,\I)-- (\LL,\I);
                    \draw (\A,\A) -- (\A,\I);
             \draw[dashed,line width=0.1mm] (\C,\A) -- (\C,\I);\draw[dashed,line width=0.1mm] (\E,\A) -- (\E,\I);\draw[dashed,line width=0.1mm] (\G,\A) -- (\G,\I);
             \draw[dashed,line width=0.1mm] (\I,\A) -- (\I,\I);\draw[dashed,line width=0.1mm] (\K,\A) -- (\K,\I);
                   \draw[solid,line width=0.4mm] (\A,\A)--(\C,\E)  -- (\G,\C) -- (\I,\G) --(\K,\C);
         \end{tikzpicture}
          \quad       \begin{tikzpicture}[scale=0.15]
            \draw (\A,\A)-- (\M,\A);
             \draw[dashed,line width=0.1mm] (\A,\E)-- (\LL,\E);
              \draw[dashed,line width=0.1mm] (\A,\C)-- (\LL,\C);
               \draw[dashed,line width=0.1mm] (\A,\G)-- (\LL,\G);
                  \draw[dashed,line width=0.1mm] (\A,\I)-- (\LL,\I);
                    \draw (\A,\A) -- (\A,\I);
             \draw[dashed,line width=0.1mm] (\C,\A) -- (\C,\I);\draw[dashed,line width=0.1mm] (\E,\A) -- (\E,\I);\draw[dashed,line width=0.1mm] (\G,\A) -- (\G,\I);
             \draw[dashed,line width=0.1mm] (\I,\A) -- (\I,\I);\draw[dashed,line width=0.1mm] (\K,\A) -- (\K,\I);
                   \draw[solid,line width=0.4mm] (\A,\A)--(\C,\E)  -- (\E,\I) -- (\I,\G) --(\K,\C);
         \end{tikzpicture}
                            \quad       \begin{tikzpicture}[scale=0.15]
            \draw (\A,\A)-- (\M,\A);
             \draw[dashed,line width=0.1mm] (\A,\E)-- (\LL,\E);
              \draw[dashed,line width=0.1mm] (\A,\C)-- (\LL,\C);
               \draw[dashed,line width=0.1mm] (\A,\G)-- (\LL,\G);
                  \draw[dashed,line width=0.1mm] (\A,\I)-- (\LL,\I);
                    \draw (\A,\A) -- (\A,\I);
             \draw[dashed,line width=0.1mm] (\C,\A) -- (\C,\I);\draw[dashed,line width=0.1mm] (\E,\A) -- (\E,\I);\draw[dashed,line width=0.1mm] (\G,\A) -- (\G,\I);
             \draw[dashed,line width=0.1mm] (\I,\A) -- (\I,\I);\draw[dashed,line width=0.1mm] (\K,\A) -- (\K,\I);
                   \draw[solid,line width=0.4mm] (\A,\A)--(\C,\E)  -- (\E,\I) -- (\G,\E) --(\K,\C);
         \end{tikzpicture}                          \quad       \begin{tikzpicture}[scale=0.15]
            \draw (\A,\A)-- (\M,\A);
             \draw[dashed,line width=0.1mm] (\A,\E)-- (\LL,\E);
              \draw[dashed,line width=0.1mm] (\A,\C)-- (\LL,\C);
               \draw[dashed,line width=0.1mm] (\A,\G)-- (\LL,\G);
                  \draw[dashed,line width=0.1mm] (\A,\I)-- (\LL,\I);
                    \draw (\A,\A) -- (\A,\I);
             \draw[dashed,line width=0.1mm] (\C,\A) -- (\C,\I);\draw[dashed,line width=0.1mm] (\E,\A) -- (\E,\I);\draw[dashed,line width=0.1mm] (\G,\A) -- (\G,\I);
             \draw[dashed,line width=0.1mm] (\I,\A) -- (\I,\I);\draw[dashed,line width=0.1mm] (\K,\A) -- (\K,\I);
                   \draw[solid,line width=0.4mm] (\A,\A)--(\C,\E)  -- (\E,\A) -- (\G,\E) --(\K,\C);
         \end{tikzpicture}
               \end{center}
         \caption{Partial knight's paths of size 5 ending at height 1.}
         \label{fig1ex1}
\end{figure}

Notice that the array $[z^nu^k](F(u)+G(u))$ corresponds to \oeis{A096587}. Figure \ref{fig1ex1} shows the partial knight's paths of size 5 ending at height 1. 

Naturally, we retrieve the generating function for the knight's paths ending on the $x$-axis $F(0)+G(0)=A(z)$.
Moreover, the generating function for the total number of partial knight's paths with respect to the size is
$$F(1)+G(1)={\frac {z(z+1)A(z)+zA_1(z)-1}{2\,{z}^{2}+2\,z-1}}.$$

The first terms of its series expansion are
$$1+z+3z^2+6z^3+16z^4+38z^5+99z^6+248z^7+646z^8+1659z^9+4342z^{10}+\ldots,$$
and the sequence of coefficients corresponds to \oeis{A096588} in \cite{Sloa}.

Since it seems difficult to obtain closed forms for $f_k$ and $g_k$, we only give the first terms of  the series expansions of $[u^k](\mathit{F(u)+G(u)})$ for the two values of $k$ ($k=1,2$), since the case $k=0$ is already mentioned above (see the series expansion of $A(x)$). For $k=1,2$ the series expansions are: 

\noindent $\bullet \quad z^2+z^3+3z^4+4z^5+12z^6+22z^7+61z^8+128z^9+\ldots,$

\noindent $\bullet \quad z+2z^3+z^4+9z^5+10z^6+42z^7+64z^8+213z^9+\ldots$

\noindent The first expansion corresponds to the sequence \oeis{A005221} in \cite{Sloa}, whereas the sequence corresponding to the second expansion is not listed in \cite{Sloa}. Notice that the cases $k=0$ and $k=1$ are already known in \cite{Lab}.

\subsection{Zigzag knight's paths}
In this part, we count partial zigzag knight's paths of a given size, i.e.,  prefixes of zigzag knight's paths ending at a given abscissa. We keep the notation of the previous section in the context of zigzag paths. 
 For $k\geq 0$, $f_k=f_k(z)$ (resp. $g_k=g_k(z)$) is the generating function where the coefficient of $z^n$ in the series expansion is the number of  partial zigzag knight's paths of size $n$  ending at height $k$ with an up-step $N$ or $E$, (resp. with a down-step  $\bar{N}$ or $\bar{E}$).

Then, we easily have the following equations:
\begin{equation}\begin{array}{l}
f_0=1,~~  f_1=z^2(f_0+g_0), \mbox{ and } f_2=z(f_0+g_0+zg_1),\\
f_k=z^2g_{k-1}+zg_{k-2}, \quad k\geq 3,\\
g_k=z^2f_{k+1}+zf_{k+2}, \quad k\geq 0.\\
\end{array}\label{eq2}
\end{equation}

Summing these recursions, we obtain:
\begin{align*}
F(u)&=1+zu(z+u)(1+G(u)),\\
G(u)&=\frac{z^2}{u}(F(u)-1) + \frac{z}{u^2}(F(u)-1-uz^2(1+g_0)),
\end{align*}
where  $g_0=G(0)$ is the generating function for the number of non-empty zigzag paths ending on the $x$-axis with respect to the size.

Solving these functional equations, we deduce
$$F(u)={\frac {g_0{u}^{2}{z}^{4}+  g_0u{z}^{5}+{u}^{2}{z}^{4}+u{z
}^{5}+{u}^{2}{z}^{3}+u{z}^{4}-{u}^{3}z-{u}^{2}{z}^{2}+u{z}^{2}+{z}^{3}
-u}{{u}^{2}{z}^{3}+u{z}^{4}+u{z}^{2}+{z}^{3}-u}}$$
and $$G(u)={\frac {{z}^{2} \left( -{u}^{2}z-u{z}^{2}+g_0z-u \right) }{u^
2z^3+uz^4+uz^2+z^3-u}}.$$
In order to compute $g_0$, we use the kernel method on $G(u)$ (see \cite{Ban2}). We have 
$$G(u)=\frac {{z}^{2} \left( -{u}^{2}z-u{z}^{2}+g_0z-u \right) }{z^3(u-s)(u-r)},$$
where 
\begin{align}\label{eqforr2}
    r={\frac {1-{z}^{4}-{z}^{2}-\sqrt {{z}^{8}-2\,{z}^{6}-{z}^{4}-2\,{z}
^{2}+1}}{2{z}^{3}}},
\end{align}
 and 
$$s={\frac {1-{z}^{4}-{z}^{2}+\sqrt {{z}^{8}-2\,{z}^{6}-{z}^{4}-2\,{
z}^{2}+1}}{2{z}^{3}}}.
$$
 It suffices to plug $u=r$ in the numerator of $G(u)$. Then, $g_0$ satisfies $  -{r}^{2}z-r{z}^{2}+g_0 z-r =0$, which implies that
 \begin{align*}g_0&=\frac{r^2z+rz^2+r}{z}=\frac{r^2z^3+rz^4+rz^2}{z^3}=\frac{r^2z^3+rz^4+rz^2+z^3}{z^3}-1\\
 &=\frac{r}{z^3}-1.
 \end{align*}
 After this, and using $sr=1$, we simplify both numerators and denominators in $F(u)$ and  $G(u)$ by extracting the common factor $(u-r)$.
 
\begin{thm}\label{theo1} The bivariate generating functions $F(u)$, $G(u)$, and $F(u)+G(u)$ are given by 
$$F(u)=\frac{ru(u+z)}{z^2(1-ru)}+1,\quad G(u)=\frac{r}{z^3(1-ru)}-1, \mbox{ and } F(u)+G(u)={\frac {r{u}^{2}z+ru{z}^{2}+r}{{z}^{3} \left( 1 -ru \right) }},$$
where $r$ is defined by Eq. (\ref{eqforr2}).
Moreover, we have
\begin{align}
f_0&=[ u^0] F(u)=1,\\
f_k&=[ u^k] F(u)=\frac{rz+1}{z^2}r^{k-1}-[k=1]\frac{1}{z^2}, ~k\geq 1,\\
g_0&=[ u^0] G(u)=\frac{r}{z^3}-1,\\
g_k&=[ u^k] G(u)=\frac{r^{k+1}}{z^3}, ~k\geq 1,
\end{align}
where $[k=1]$ equals $1$ whenever $k=1$ and $0$ otherwise. 
\end{thm}

The first terms of the series expansion of $F(u) + G(u)$ are
\begin{multline*}
1 + u^2 z + (1 + u) z^2 + (u + u^2) z^3 + (2 + u + u^3) z^4 + (2 u + 
    3 u^2) z^5 \\+ (4 + 2 u + u^2 + 2 u^3) z^6 + (5 u + 6 u^2 + u^4) z^7\ldots
\end{multline*}

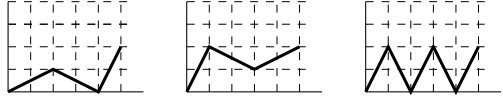
\begin{figure}[h]
 \begin{center}
        \begin{tikzpicture}[scale=0.15]
            \draw (\A,\A)-- (\M,\A);
             \draw[dashed,line width=0.1mm] (\A,\E)-- (\LL,\E);
              \draw[dashed,line width=0.1mm] (\A,\C)-- (\LL,\C);
               \draw[dashed,line width=0.1mm] (\A,\G)-- (\LL,\G);
            \draw[dashed,line width=0.1mm] (\A,\G)-- (\LL,\G);
             \draw[dashed,line width=0.1mm] (\A,\I)-- (\LL,\I);
                    \draw (\A,\A) -- (\A,\I);
             \draw[dashed,line width=0.1mm] (\C,\A) -- (\C,\I);\draw[dashed,line width=0.1mm] (\E,\A) -- (\E,\I);\draw[dashed,line width=0.1mm] (\G,\A) -- (\G,\I);
             \draw[dashed,line width=0.1mm] (\I,\A) -- (\I,\I);\draw[dashed,line width=0.1mm] (\K,\A) -- (\K,\I);
                   \draw[solid,line width=0.4mm] (\A,\A)--(\E,\C)  -- (\I,\A) -- (\K,\E);
         \end{tikzpicture}
          \quad       \begin{tikzpicture}[scale=0.15]
            \draw (\A,\A)-- (\M,\A);
             \draw[dashed,line width=0.1mm] (\A,\E)-- (\LL,\E);
              \draw[dashed,line width=0.1mm] (\A,\C)-- (\LL,\C);
               \draw[dashed,line width=0.1mm] (\A,\G)-- (\LL,\G);
                  \draw[dashed,line width=0.1mm] (\A,\I)-- (\LL,\I);
                    \draw (\A,\A) -- (\A,\I);
             \draw[dashed,line width=0.1mm] (\C,\A) -- (\C,\I);\draw[dashed,line width=0.1mm] (\E,\A) -- (\E,\I);\draw[dashed,line width=0.1mm] (\G,\A) -- (\G,\I);
             \draw[dashed,line width=0.1mm] (\I,\A) -- (\I,\I);\draw[dashed,line width=0.1mm] (\K,\A) -- (\K,\I);
                   \draw[solid,line width=0.4mm] (\A,\A)--(\C,\E)  -- (\G,\C) -- (\K,\E) ;
         \end{tikzpicture}
                            \quad       \begin{tikzpicture}[scale=0.15]
            \draw (\A,\A)-- (\M,\A);
             \draw[dashed,line width=0.1mm] (\A,\E)-- (\LL,\E);
              \draw[dashed,line width=0.1mm] (\A,\C)-- (\LL,\C);
               \draw[dashed,line width=0.1mm] (\A,\G)-- (\LL,\G);
                  \draw[dashed,line width=0.1mm] (\A,\I)-- (\LL,\I);
                    \draw (\A,\A) -- (\A,\I);
             \draw[dashed,line width=0.1mm] (\C,\A) -- (\C,\I);\draw[dashed,line width=0.1mm] (\E,\A) -- (\E,\I);\draw[dashed,line width=0.1mm] (\G,\A) -- (\G,\I);
             \draw[dashed,line width=0.1mm] (\I,\A) -- (\I,\I);\draw[dashed,line width=0.1mm] (\K,\A) -- (\K,\I);
                   \draw[solid,line width=0.4mm] (\A,\A)--(\C,\E)  -- (\E,\A) -- (\G,\E) -- (\I, \A)-- (\K, \E) ;
         \end{tikzpicture}                      
               \end{center}
         \caption{Partial zigzag knight's paths of size 5 ending at height 2.}
         \label{fig1ex2}
\end{figure}

Figure \ref{fig1ex2} shows the partial zigzag knight's paths of length 5 ending at height 2. 

Here are other examples of  the first terms of the series expansions of $[ u^k](\mathit{F(u)+G(u)})$ for $k=1,2,3$:

\noindent $\bullet \quad z^2+z^3+z^4+2z^5+2z^6+5z^7+4z^8+12z^9+8z^{10}+28z^{11}+17z^{12}+\ldots,$

\noindent $\bullet \quad z+z^3+3z^5+z^6+6z^7+3z^8+13z^9+9z^{10}+29z^{11}+25z^{12}+65z^{13}+\ldots,$

\noindent $\bullet \quad z^4+2z^6+6z^8+z^9+15z^{10}+4z^{11}+37z^{12}+14z^{13}+91z^{14}+44z^{15}+\ldots,$

\noindent which do not appear in \cite{Sloa}.

\begin{cor} The generating function for the number of zigzag knight's paths of a given size ending on the $x$-axis is 
$$F(0)+G(0)=\frac{r}{z^3},$$ and the sequence of non-zero coefficients in the series expansion corresponds to the generalized Catalan numbers (see \oeis{A004148} in \cite{Sloa}).  The first terms of the series expansion are $$1+z^2+2z^4+4z^6+8z^8+17z^{10}+37z^{12}+82z^{14}+185z^{16}+423z^{18}+978z^{20}+\ldots$$\label{cor1}
\end{cor}

\begin{cor} The generating function for the number of partial zigzag knight's paths of a given size is
$$F(1)+G(1)={\frac {r{z}^{2}+rz+r}{{z}^{3} \left( 1-r \right) }}, 
$$ and the sequence of coefficients of $z^n$ in the series expansion corresponds to \oeis{A088518} in \cite{Sloa}, which  counts also symmetric secondary structures of RNA molecules with $n$ nucleotides. The first terms of the series expansions are  $$1+z+2z^2+2z^3+4z^4+5z^5+9z^6+12z^7+21z^8+29z^{9}+50z^{10}+\ldots$$
	\end{cor}
 The next theorem provides close forms for the coefficients of $z^{n}$ in the series expansions of $f_k$ and $g_k$, $k\geq 0$.

\begin{thm} We have $[z^0]f_0=1$, and for $n\geq 1$, $k\geq 1$, 
\begin{align*}
[z^{2n}]f_{2k+1}&=a(n-k,2k+1)+a(n-k+1,2k),\\ 
[z^{2n-1}]f_{2k}&=a(n-k,2k)+a(n-k+1,2k-1),\\
[z^{2n}]g_{2k}&=a(n-k+1,2k+1),\\
[z^{2n-1}]g_{2k-1}&=a(n-k+1,2k),
\end{align*}
where $$a(n,k)=\sum_{i=0}^{\lfloor\frac{n-1}{2} \rfloor}\frac{k}{n-i}\binom{n-i}{i+ k}\binom{n-i}{i}.
$$
All other coefficients are equal to $0$.
\end{thm}
\noindent {\it Proof.} 
Let $\hat{r}$ be the function $z^{-1/2}r(z^{1/2})$, where $r(z)$ is defined in  \eqref{eqforr2}. Then
$$\hat{r}=\frac{1-z-z^2-\sqrt{1 - 2 z - z^2 - 2 z^3 + z^4}}{2z^2}.$$
The function $\hat{r}$ satisfies the functional equation 
 $\hat{r}=z(1+(1+z)\hat{r}+z\hat{r}^2)$. 
 Consider the auxiliary function $q(z,t)$ defined by
 $$q(z,t)=z(1+(1+t)q(z,t)+tq(z,t)^2)=z\Phi(q(z,t)),$$
 where $\Phi(u)=(1+u)(1+tu)$.  From  Lagrange inversion, (see \cite{Merl} for instance), we have that 
\begin{align*}
[z^n]q(z,t)^k&=\frac{k}{n}[u^{n-k}]\Phi(u)^n=\frac{k}{n}[u^{n-k}](u+1)^{n}(tu+1)^n\\
&=\frac{k}{n}\sum_{i=0}^{n-k}\binom{n}{n-k-i}\binom{n}{i}t^i, \quad n\geq 1.
\end{align*} 
So, we obtain
\begin{align*}
[z^n]\hat{r}^k&=[z^n]q(z,z)^k=[z^n]\sum_{\ell\geq 1}\frac{k}{\ell}\sum_{i=0}^{\ell-k}\binom{\ell}{\ell-k-i}\binom{\ell}{i}z^{\ell+i}\\
&=[z^n]\sum_{i\geq 0}\sum_{\ell\geq 0}\frac{k}{\ell+i+1}\binom{\ell + i + 1}{\ell + 1 - k}\binom{\ell + i + 1}{i}z^{\ell+2i+1}.
\end{align*} 
Setting $h = \ell + 2i + 1$, this implies
\begin{align*}
[z^n]\hat{r}^k&=[z^n]\sum_{i\geq 0}\sum_{h\geq 2i+1}\frac{k}{h-i}\binom{h-i}{h - 2 i - k}\binom{h-i}{i}z^{h}\\
&=\sum_{i=0}^{\lfloor\frac{n-1}{2} \rfloor}\frac{k}{n-i}\binom{n-i}{ i + k}\binom{n-i}{i}=a(n,k).
\end{align*} 
Then, for all positive integers $n $ and $k$, we have $$[z^{2n}]r^{2k}=a(n-k,2k)\quad \text{ and } \quad [z^{2n-1}]r^{2k-1}=a(n-k,2k-1).$$ 
From  Theorem \ref{theo1}, we conclude
\begin{align*}
[z^{2n}]f_{2k+1}&=[z^{2n}]\left(\frac{r^{2k+1}}{z} +\frac{r^{2k}}{z^2}\right)=[z^{2n+1}]r^{2k+1} +[z^{2n+2}]r^{2k}\\
&=a(n-k,2k+1)+a(n-k+1,2k).
\end{align*}
Analogously, we obtain the remaining identities.
\hfill $\Box$

\subsection{A bijective approach}

Corollary~\ref{cor1} proves that the set of zigzag knight's paths of size $2n$ ending on the $x$-axis is equinumerous to the set of peakless Motzkin paths of length $n+1$, i.e., lattice paths in $\Bbb{N}^2$, starting at the origin, ending at $(n+1,0)$, consisting  of steps $U=(1,1)$, $D=(1,-1)$, $F=(1,0)$, and avoiding the pattern $UD$ (see for instance \cite{Cam}). Any non-empty peakless Motzkin path is either of the form (1) $FQ$ or (2) $URDS$ or (3) $URD$, where $Q$ (resp. $R$, $S$) is a possible empty (resp. non-empty) peakless Motzkin path. Let $\mathcal{M}_n$ be the set of peakless Motzkin paths of length $n$, and $\mathcal{M}=\bigcup_{n\geq 0}\mathcal{M}_n$.

\begin{figure}[H]
\centering
\begin{tikzpicture}[scale=0.26]
\draw (1,0) node [gray,scale=1.6]{$\epsilon$};

\draw[very thick] (0.5,-10) -- (1.5,-10);

\draw[orange, very thick] (11,0) .. controls (12,3) and (16,3) .. (17,0) node[black,xshift=-24,yshift=6] {\large $\beta$};
\draw[very thick] (7,0) -- (9,1);
\draw[very thick] (9,1) -- (11,0);

\draw[orange, very thick] (8,-9) .. controls (9,-6) and (13,-6) .. (14,-9) node[black,xshift=-24,yshift=8] {\large $\psi(\beta)$};
\draw[very thick] (7,-10) -- (8,-9);
\draw[very thick] (14,-9) -- (15,-10);
\draw[gray, thin, dashed] (8,-9) -- (14,-9);

\draw[orange, very thick] (22,0) .. controls (23,3) and (27,3) .. (28,0) node[black,xshift=-24,yshift=6] {\large $\beta$};
\draw[very thick] (20,0) -- (21,2);
\draw[very thick] (21,2) -- (22,0);

\draw[orange, very thick] (21,-10) .. controls (22,-7) and (26,-7) .. (27,-10) node[black,xshift=-24,yshift=6] {\large $\psi(\beta)$};
\draw[very thick] (20,-10) -- (21,-10);
 
\draw[very thick] (30,0) -- (31,2);
\draw[very thick] (31,2) -- (33,1);
\draw[orange, very thick] (33,1) .. controls (34,4) and (38,4) .. (39,1) node[black,xshift=-24,yshift=6] {\large $\beta$};
\draw[very thick] (39,1) -- (41,2);
\draw[very thick] (41,2) -- (42,0);
\draw[orange, very thick] (42,0) .. controls (43,3) and (47,3) .. (48,0) node[black,xshift=-24,yshift=6] {\large $\gamma$};

\draw[orange, very thick] (32.5,-9) .. controls (33.5,-6) and (37.5,-6) .. (38.5,-9) node[black,xshift=-24,yshift=6] {\large $\psi(\beta)$};

\draw[orange, very thick] (39.5,-10) .. controls (40.5,-7) and (44.5,-7) .. (45.5,-10) node[black,xshift=-24,yshift=6] {\large $\psi(\gamma)$};
\draw[very thick] (31.5,-10) -- (32.5,-9);
\draw[very thick] (38.5,-9) -- (39.5,-10);
\draw[gray, thin, dashed] (32.5,-9) -- (38.5,-9);

\draw[very thick,|->] (1,-1.7) -- (1,-5.5);
\draw[very thick,|->] (11,-1.7) -- (11,-5.5);
\draw[very thick,|->] (23.5,-1.7) -- (23.5,-5.5);
\draw[very thick,|->] (38.5,-1.7) -- (38.5,-5.5);
\end{tikzpicture}
\caption{Illustration of the map $\psi$ according to Definition~\ref{Def1}.}
\label{fig3}
\end{figure}
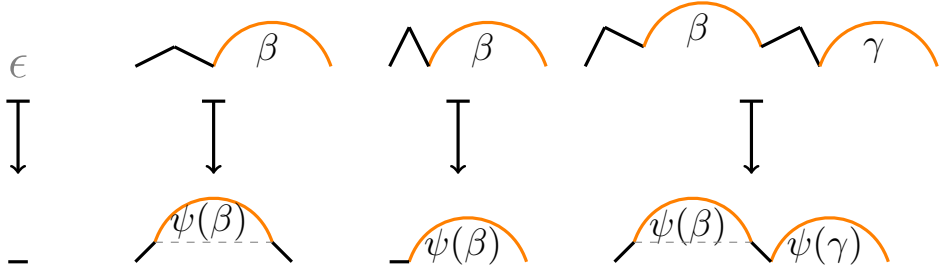

 On the other hand, any non-empty zigzag knight's path is either of the form (1) $N\bar{N}\beta$, (2) $E\bar{E}\beta$ or  (3) $N\bar{E} \beta E\bar{N}\gamma$, where $\beta,\gamma$ are possibly empty zigzag knight's paths, and $N=(1,2)$, $\bar{N}=(1,-2)$, $E=(2,1)$, $\bar{E}=(2,-1)$.  We refer to the top of  Figure \ref{fig3} for an illustration of these different forms.  Let $\mathcal{R}_n$ be the set of zigzag paths  of size $2n$ ending on the $x$-axis, and $\mathcal{R}=\bigcup_{n\geq 0}\mathcal{R}_n$.

\begin{defn} We  recursively define the map $\psi$ from $\mathcal{R}$ to $\mathcal{M}$ as follows. For $\alpha\in\mathcal{R}$, we set:
$$
\psi(\alpha)=\left\{\begin{array}{llr}
F&\text{if }\alpha=\epsilon,&(i)\\
U\psi(\beta)D&\text{if }\alpha=E\bar{E}\beta \text{ with } \beta\in\mathcal{R},&(ii)\\
F\psi(\beta)&\text{if }\alpha=N\bar{N}\beta \text{ with } \beta\in\mathcal{R},&(iii)\\
U\psi(\beta)D\psi(\gamma)&\text{if }\alpha=N\bar{E}\beta E\bar{N}\gamma \text{ with } \beta,\gamma \in\mathcal{R}. &(iv)
\end{array}\right.
$$
\label{Def1}
\end{defn}
Due to the recursive definition, the image  of a zigzag knight's path of 
size $2n$ under $\psi$ is a peakless Motzkin path of length $n+1$. For instance, the images of $N\bar{N}$, $E\bar{E}$ and $N\bar{N}N\bar{N}$ are $FF$, $UFD$ and $FFF$, respectively. We refer to Figure~\ref{fig3} for an illustration of this mapping.

\begin{thm}
For all $n\geq 0$, the map $\psi$ induces a bijection between $\mathcal{R}_n$ and $\mathcal{M}_{n+1}$.
\label{bij}
\end{thm}
\noindent {\it Proof.} Due to Corollary~1 in Section 2.2, the cardinality of $\mathcal{R}_n$ is given by the $n$-th term of generalized Catalan number. So it suffices to prove the injectivity of $\psi$. We proceed by induction on $n$. The case $n=0$ is obvious since $\mathcal{R}_0=\{\epsilon\}$ and $\mathcal{M}_1=\{F\}$.
 For all $k\leq n$, we assume that $\psi$ is an injection from $\mathcal{R}_k$ to $\mathcal{M}_{k+1}$, and we prove the result for $k=n+1$.
 According to Definition~\ref{Def1}, if $\alpha$ and $\beta$ in $\mathcal{R}_{n+1}$ satisfy $\psi(\alpha)=\psi(\beta)$, then $\alpha$ and $\beta$ necessarily come from the same case among ($i$)~--~($iv$). Using the induction hypothesis, we conclude directly that $\alpha=\beta$, which completes the induction. Thus $\psi$ is injective and so bijective. \hfill $\Box$
\medskip

For instance, the image of $E\bar{E}N\bar{N}N\bar{E}N\bar{N}E\bar{N}$ under $\psi$ is $UFUFFDFD$.

\section{Knight's paths of a given length}

Knight's paths of a given length $n$ are clearly in one-to-one correspondence with basketball walks (see \cite{Ayy,Bank, Bett, Bousq}), i.e., lattices paths in $\Bbb{N}^2$, starting at the origin, and consisting of $n$ steps lying in $S'=\{(1,1),(1,-1), (1,2),(1,-2)\}$. The correspondence is obtained by replacing $(1,1)$ and $(1,-1)$ with $E$ and $\bar{E}$, respectively.   

In \cite{Bank}, Banderier et al. (see Prop. 3.7) prove that the generating function $E(z)=\sum_{n\geq 0}e_nz^n$ for the number of these paths with respect to the length satisfies 
$$z^4E^4-(2z^3+z^2)E^3+(3z^2+2z)E^2-(2z+1)E+1=0,$$
and they obtain
$$e_n=\frac{1}{n+1}\sum\limits_{i=0}^{\lfloor n/2\rfloor}{2n+2\choose i}{n-i-1\choose n-2i}.$$ The first values of $e_n$, $n\geq 0$, are $1,0,2,2,11,24,93,272,971,3194,11293$ (sequence \oeis{A187430}).

\subsection{Partial knight's paths}

Partial knight's paths are also handled by Banderier et al. in \cite{Bank}. In their Prop. 3.2, they prove that the generating function for basketball walks with respect to the number of steps, starting at the origin, ending at height $k$, is given by $$\frac{u_1^{k+1}(z)-u_2^{k+1}(z)}{u_1(z)-u_2(z)},$$
where $u_1(z)$ and $u_2(z)$ are the small roots of $$1-z(u^{-2}+u^{-1}+u+u^2)=0.$$

For instance (see Eq. (3.11) and (3.12) in \cite{Bank}), the generating functions for basketball walks with respect to the number of steps, starting at the origin, ending at height $k$, $k=1,2$, are respectively 

$$\frac{1}{2}\left(-1+\sqrt{\frac{2-3z-2\sqrt{1-4z}}{z}}\right)\mbox{ \ and \ } \frac{3-\sqrt{1-4z}-\sqrt{2+12z+2\sqrt{1-4z}}}{4z}.$$

The first terms of the series expansions are respectively 

\noindent $\bullet \quad z+z^2+3z^3+7z^4+22z^5+65z^6+213z^7+\ldots$ (sequence \oeis{A166135}) and 

\noindent $\bullet \quad z+z^2+4z^3+9z^4+31z^5+91z^6+309z^7+\ldots$ (sequence \oeis{A111160}).

\subsection{Zigzag knight's paths}

In this part, we count partial zigzag knight's paths of a given length. We adopt the notation used in Section 2.2 for  zigzag paths of given length.  For $k\geq 0$, $f_k=f_k(z)$ (resp. $g_k=g_k(z)$) is the generating function where the coefficient of $z^n$ in the series expansion is the number of  partial zigzag knight's paths of length $n$. Then, we easily have the following equations:
\begin{equation}\begin{array}{l}
f_0=1,  f_1=z(f_0+g_0), \mbox{ and } f_2=z(f_0+g_0+g_1),\\
f_k=zg_{k-1}+zg_{k-2}, \quad k\geq 3,\\
g_k=zf_{k+1}+zf_{k+2}, \quad k\geq 0.
\end{array}\label{equ1}
\end{equation}

Summing the recursions in (\ref{equ1}), we have:
\begin{align*}
F(u)&=1+zu(1+u)(1+G(u)),\\
G(u)&=\frac{z}{u}(F(u)-1) + \frac{z}{u^2}(F(u)-1-uz(1+g_0)),
\end{align*}
where  $g_0=G(0)$ is the generating function with respect to the length for the number of non-empty zigzag knight's paths ending on the $x$-axis.

Solving these functional equations, we deduce
$$F(u)={\frac {g_0{u}^{2}{z}^{3} + g_0u{z}^{3}+{u}^{2}{z}^{3}-{u}
^{3}z+{u}^{2}{z}^{2}+u{z}^{3}-{u}^{2}z+2\,u{z}^{2}+{z}^{2}-u}{{u}^{2}{
z}^{2}+2\,u{z}^{2}+{z}^{2}-u}}
,$$
and $$G(u)={\frac {{z}^{2} \left( -{u}^{2} + g_0 -2 u \right) }{{u}^{2}{z}^{2}
+2\,u{z}^{2}+{z}^{2}-u}}.$$
In order to compute $g_0$, we use the kernel method on $G(u)$. We have 
$$G(u)=\frac{{z}^{2} \left( -{u}^{2} + g_0-2\,u \right)}{z^2(u-r)(u-s)},$$
where 
\begin{align}\label{eqforr}
r={\frac {1-2\,{z}^{2}-\sqrt {1-4\,{z}^{2}}}{2{z}^{2}}}
\end{align}
and $$s={\frac {1-2\,{z}^{2}+\sqrt {1-4\,{z}^{2}}}{2{z}^{2}}}.$$

Plugging $u=r$ in the numerator of $G(u)$, we obtain $-r^2 + g_0 - 2r=0$, which implies that $g_0=r(2+r)$. After this, and using $sr=1$, we simplify both numerators and denominators in $F(u)$ and $G(u)$ by extracting the common factor $(u-r)$.
\begin{thm} The bivariate generating functions $F(u)$, $G(u)$, and $F(u)+G(u)$ are given by 
$$F(u)={\frac {r{u}^{2}+ur \left( 1-z \right) +z}{ z\left( 1-ru \right)}}
, \quad G(u)=\frac{r(u+r+2)}{1-ru},
$$
and
$$F(u)+G(u)={\frac {{r}^{2}z+r{u}^{2}+ru+2\,rz+z}{ z\left( 1-ru \right)}}
,$$
where $r$ is defined by Eq. (\ref{eqforr}).
Moreover, we have 
\begin{align}
f_0&=[u^0]F(u)=1,\\
f_k&=[u^k]F(u)=\frac{1+r}{z}r^{k-1}
-[k=1]\frac{1}{z},~k\geq 1, \label{eq10fk}\\
g_0&=[u^0]G(u)=r(2+r),\\
g_k&=[u^k]G(u)=\frac{r^{k+1}}{z^2}, ~k\geq 1.
\end{align}
\end{thm}

The first terms of the series expansion of $F(u) + G(u)$ are
\begin{multline*}
1 + (u + u^2) z + (2 + u) z^2 + (2 u + 3 u^2 + u^3) z^3 + (5 + 4 u + 
    u^2) z^4 \\ + (5 u + 9 u^2 + 5 u^3 + u^4) z^5 + (14 + 14 u + 6 u^2 + 
    u^3) z^6+\ldots
\end{multline*}

\begin{figure}[h]
 \begin{center}
        \begin{tikzpicture}[scale=0.15]
            \draw (\A,\A)-- (\W,\A);
             \draw[dashed,line width=0.1mm] (\A,\E)-- (\V,\E);
              \draw[dashed,line width=0.1mm] (\A,\C)-- (\V,\C);
               \draw[dashed,line width=0.1mm] (\A,\G)-- (\V,\G);
            \draw[dashed,line width=0.1mm] (\A,\G)-- (\V,\G);
             \draw[dashed,line width=0.1mm] (\A,\I)-- (\V,\I);
                    \draw (\A,\A) -- (\A,\I);
             \draw[dashed,line width=0.1mm] (\C,\A) -- (\C,\I);\draw[dashed,line width=0.1mm] (\E,\A) -- (\E,\I);\draw[dashed,line width=0.1mm] (\G,\A) -- (\G,\I);
             \draw[dashed,line width=0.1mm] (\I,\A) -- (\I,\I);
             \draw[dashed,line width=0.1mm] (\K,\A) -- (\K,\I);
             \draw[dashed,line width=0.1mm] (\M,\A) -- (\M,\I);
                \draw[dashed,line width=0.1mm] (\O,\A) -- (\O,\I);
            \draw[dashed,line width=0.1mm] (\Q,\A) -- (\Q,\I); 
            \draw[dashed,line width=0.1mm] (\S,\A) -- (\S,\I);
            \draw[dashed,line width=0.1mm] (\U,\A) -- (\U,\I);
                   \draw[solid,line width=0.4mm] (\A,\A)--(\E,\C)  -- (\I,\A) -- (\M,\C)-- (\Q,\A)-- (\U,\C);
         \end{tikzpicture}
          \quad       \begin{tikzpicture}[scale=0.15]
            \draw (\A,\A)-- (\S,\A);
             \draw[dashed,line width=0.1mm] (\A,\E)-- (\R,\E);
              \draw[dashed,line width=0.1mm] (\A,\C)-- (\R,\C);
               \draw[dashed,line width=0.1mm] (\A,\G)-- (\R,\G);
                  \draw[dashed,line width=0.1mm] (\A,\I)-- (\R,\I);
                    \draw (\A,\A) -- (\A,\I);
             \draw[dashed,line width=0.1mm] (\C,\A) -- (\C,\I);\draw[dashed,line width=0.1mm] (\E,\A) -- (\E,\I);\draw[dashed,line width=0.1mm] (\G,\A) -- (\G,\I);
             \draw[dashed,line width=0.1mm] (\I,\A) -- (\I,\I);\draw[dashed,line width=0.1mm] (\K,\A) -- (\K,\I);
               \draw[dashed,line width=0.1mm] (\M,\A) -- (\M,\I);
                \draw[dashed,line width=0.1mm] (\O,\A) -- (\O,\I);
            \draw[dashed,line width=0.1mm] (\Q,\A) -- (\Q,\I); 
                    \draw[solid,line width=0.4mm] (\A,\A)--(\E,\C)  -- (\I,\A) -- (\K,\E)-- (\M,\A)-- (\Q,\C);
         \end{tikzpicture}
                            \quad     
         \begin{tikzpicture}[scale=0.15]
            \draw (\A,\A)-- (\S,\A);
             \draw[dashed,line width=0.1mm] (\A,\E)-- (\R,\E);
              \draw[dashed,line width=0.1mm] (\A,\C)-- (\R,\C);
               \draw[dashed,line width=0.1mm] (\A,\G)-- (\R,\G);
                  \draw[dashed,line width=0.1mm] (\A,\I)-- (\R,\I);
                    \draw (\A,\A) -- (\A,\I);
             \draw[dashed,line width=0.1mm] (\C,\A) -- (\C,\I);\draw[dashed,line width=0.1mm] (\E,\A) -- (\E,\I);\draw[dashed,line width=0.1mm] (\G,\A) -- (\G,\I);
             \draw[dashed,line width=0.1mm] (\I,\A) -- (\I,\I);\draw[dashed,line width=0.1mm] (\K,\A) -- (\K,\I);
               \draw[dashed,line width=0.1mm] (\M,\A) -- (\M,\I);
                \draw[dashed,line width=0.1mm] (\O,\A) -- (\O,\I);
            \draw[dashed,line width=0.1mm] (\Q,\A) -- (\Q,\I); 
                    \draw[solid,line width=0.4mm] (\A,\A)--(\C,\E)  -- (\G,\C) -- (\K,\E)-- (\M,\A)-- (\Q,\C);
         \end{tikzpicture}
                            \quad  
           \begin{tikzpicture}[scale=0.15]
            \draw (\A,\A)-- (\S,\A);
             \draw[dashed,line width=0.1mm] (\A,\E)-- (\R,\E);
              \draw[dashed,line width=0.1mm] (\A,\C)-- (\R,\C);
               \draw[dashed,line width=0.1mm] (\A,\G)-- (\R,\G);
                  \draw[dashed,line width=0.1mm] (\A,\I)-- (\R,\I);
                    \draw (\A,\A) -- (\A,\I);
             \draw[dashed,line width=0.1mm] (\C,\A) -- (\C,\I);\draw[dashed,line width=0.1mm] (\E,\A) -- (\E,\I);\draw[dashed,line width=0.1mm] (\G,\A) -- (\G,\I);
             \draw[dashed,line width=0.1mm] (\I,\A) -- (\I,\I);\draw[dashed,line width=0.1mm] (\K,\A) -- (\K,\I);
               \draw[dashed,line width=0.1mm] (\M,\A) -- (\M,\I);
                \draw[dashed,line width=0.1mm] (\O,\A) -- (\O,\I);
            \draw[dashed,line width=0.1mm] (\Q,\A) -- (\Q,\I); 
                    \draw[solid,line width=0.4mm] (\A,\A)--(\C,\E)  -- (\E,\A) -- (\I,\C)-- (\M,\A)-- (\Q,\C);
         \end{tikzpicture}                            \quad  
                            \begin{tikzpicture}[scale=0.15]
            \draw (\A,\A)-- (\O,\A);
             \draw[dashed,line width=0.1mm] (\A,\E)-- (\N,\E);
              \draw[dashed,line width=0.1mm] (\A,\C)-- (\N,\C);
               \draw[dashed,line width=0.1mm] (\A,\G)-- (\N,\G);
                  \draw[dashed,line width=0.1mm] (\A,\I)-- (\N,\I);
                                          \draw (\A,\A) -- (\A,\I);
         \draw[dashed,line width=0.1mm] (\C,\A) -- (\C,\I);\draw[dashed,line width=0.1mm] (\E,\A) -- (\E,\I);\draw[dashed,line width=0.1mm] (\G,\A) -- (\G,\I);
             \draw[dashed,line width=0.1mm] (\I,\A) -- (\I,\I);\draw[dashed,line width=0.1mm] (\K,\A) -- (\K,\I);
             \draw[dashed,line width=0.1mm] (\M,\A) -- (\M,\I);
                             \draw[solid,line width=0.4mm] (\A,\A)--(\C,\E)  -- (\E,\A) -- (\G,\E) -- (\I, \A)-- (\M, \C) ;
         \end{tikzpicture}                      
               \end{center}
         \caption{Partial zigzag knight's path of length 5 ending at height 1.}
         \label{fig1ex3}
\end{figure}
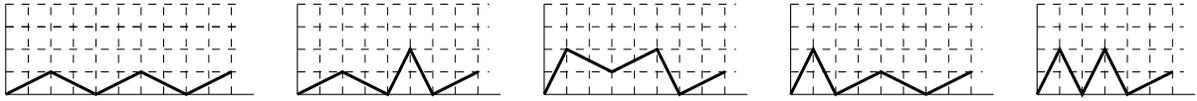

Figure \ref{fig1ex3} shows the partial zigzag knight's paths of length 5 ending at height 1. 

\begin{cor} The generating function for the number of zigzag knight's paths of a given length ending on the $x$-axis is 
$$F(0)+G(0)=1+r(2+r),$$
and the sequence of non-zero coefficients in the series expansion corresponds to the Catalan numbers (\oeis{A000108} in \cite{Sloa}). The first terms of the series expansion are $$1+2z^2+5z^4+14z^6+42z^{8}+132z^{10}+429z^{12}+1430z^{14}+4862z^{16}+16796z^{18}+58786x^{20}+\ldots$$
\label{cor3}
\end{cor}

\begin{cor} The generating function for the number of partial zigzag knight's paths of a given length is
$$F(1)+G(1)={\frac {2\,rz+r}{ \left( 1-r \right) {z}^{2}}}
={\frac {4\,z+2}{1-4\,{z}^{2}+\sqrt {1-4\,{z}^{2}}}},$$
and the sequence of coefficients in the series expansion corresponds to a shift of \oeis{A001405} in \cite{Sloa}. The first terms of the series expansion are $$1+2z+3z^2+6z^3+10z^4+20z^5+35z^6+70z^7+126z^8+252z^9+462z^{10}+\ldots$$
\end{cor}

\begin{cor} We have $[z^0]f_0=1$ and, for $n\geq 1$,
\begin{align*}
[z^{2n-1}]f_k&=\frac{2k-1}{n+k}{2n \choose n-k+1},~k\geq 1,\\
[z^{2n}]g_0&=\frac{1}{n+2}{2n+2 \choose n+1},\\
[z^{2n}]g_k&=\frac{k+1}{n+1}\binom{2n+2}{n-k}, ~k\geq 1.
\end{align*}
All other coefficients are equal to $0$.
\end{cor}
\noindent {\it Proof.} 
The generating function $\hat{r}:=r(z^{1/2})$ (see Eq. \eqref{eqforr}) satisfies the functional equation $\hat{r}=z\Phi(\hat{r}),$ where $\Phi(u)=u^2+2u+1$. From the Lagrange inversion we have that 
\begin{align*}
[z^{2n}]r^k=[z^n]\hat{r}^k=\frac{k}{n}[u^{n-k}]\Phi(u)^n=\frac{k}{n}[u^{n-k}](u+1)^{2n}=\frac{k}{n}\binom{2n}{n-k}, \quad n\geq 1.
\end{align*} 
From Eq. \eqref{eq10fk} we have
\begin{align*}
[z^{2n-1}]f_k&=[z^{2n}](r^{k-1}+r^{k})= \frac{k-1}{n}\binom{2n}{n-k+1}+\frac{k}{n}\binom{2n}{n-k}\\
&=\frac{k-1}{n}\binom{2n}{n-k+1}+\frac{k(n-k+1)}{n(n+k)}\binom{2n}{n-k+1}\\
&=\frac{2k-1}{n+k}{2n \choose n-k+1}.
\end{align*}
On the other hand,
$$[z^{2n}]g_0=[z^{2n}]r(2+r)=\frac{2}{n}\binom{2n}{n-1} + \frac{2}{n}\binom{2n}{n-2}=\frac{1}{n+2}{2n+2 \choose n+1},$$
and, for $k\geq 1$
$$[z^{2n}]g_k=[z^{2n}]\frac{r^{k+1}}{z^2}=[z^{2n+2}]r^{k+1}=\frac{k+1}{n+1}\binom{2n+2}{n-k}.$$
\hfill $\Box$

Here are other examples of  the first terms of the series expansions of $[u^k](\mathit{F(u)+G(u)})$ for $k=1,2,3$:

\noindent $\bullet \quad z+z^2+2z^3+4z^4+5z^5+14z^6+14z^7+48z^8+42z^9+165z^{10}+\ldots,$

\noindent $\bullet \quad z+3z^3+z^4+9z^5+6z^6+28z^7+27z^8+90z^9+110z^{10}+\ldots,$

\noindent $\bullet \quad z^3+5z^5+z^6+20z^7+8z^8+75z^9+44z^{10}+275z^{11}+208z^{12}+\ldots,$

\noindent which do not appear in \cite{Sloa}.

\subsection{A bijective approach}

Corollary~\ref{cor3} proves that there is a bijection between the set of zigzag knight's paths of length $2n$ ending on the $x$-axis and the set of Dyck paths of  length $2(n+1)$. Any non-empty Dyck path is either of the form (1) $UD$, (2) $UQD$, (3) $UDQ$ or (4) $UQDR$, where $Q$ and $R$ are non-empty Dyck paths. Let $\mathcal{D}_n$ (resp. $\mathcal{S}_n$) be the set of Dyck paths (resp. zigzag paths) of length $2n$, and $\mathcal{D}=\bigcup_{n\geq 0}\mathcal{D}_n$ (resp. $\mathcal{S}=\bigcup_{n\geq 0}\mathcal{S}_n$).

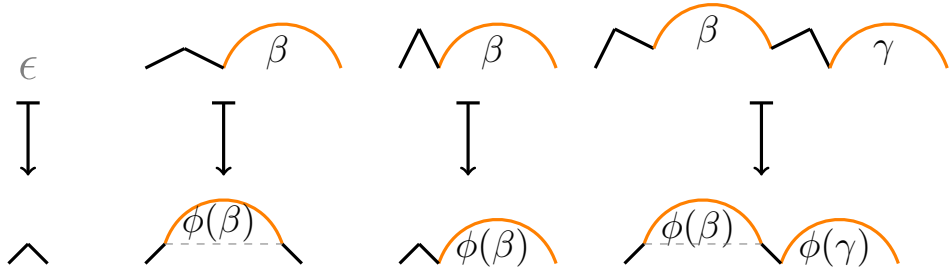
\begin{figure}[H]
\centering
\begin{tikzpicture}[scale=0.26]
\draw (1,0) node [gray,scale=1.6]{$\epsilon$};

\draw[very thick] (0,-10) -- (1,-9);\draw[very thick] (1,-9) -- (2,-10);

\draw[orange, very thick] (11,0) .. controls (12,3) and (16,3) .. (17,0) node[black,xshift=-24,yshift=6] {\large $\beta$};
\draw[very thick] (7,0) -- (9,1);
\draw[very thick] (9,1) -- (11,0);

\draw[orange, very thick] (8,-9) .. controls (9,-6) and (13,-6) .. (14,-9) node[black,xshift=-24,yshift=8] {\large $\phi(\beta)$};
\draw[very thick] (7,-10) -- (8,-9);
\draw[very thick] (14,-9) -- (15,-10);
\draw[gray, thin, dashed] (8,-9) -- (14,-9);

\draw[orange, very thick] (22,0) .. controls (23,3) and (27,3) .. (28,0) node[black,xshift=-24,yshift=6] {\large $\beta$};
\draw[very thick] (20,0) -- (21,2);
\draw[very thick] (21,2) -- (22,0);

\draw[orange, very thick] (22,-10) .. controls (23,-7) and (27,-7) .. (28,-10) node[black,xshift=-24,yshift=6] {\large $\phi(\beta)$};
\draw[very thick] (20,-10) -- (21,-9);\draw[very thick] (21,-9) -- (22,-10);
 
\draw[very thick] (30,0) -- (31,2);
\draw[very thick] (31,2) -- (33,1);
\draw[orange, very thick] (33,1) .. controls (34,4) and (38,4) .. (39,1) node[black,xshift=-24,yshift=6] {\large $\beta$};
\draw[very thick] (39,1) -- (41,2);
\draw[very thick] (41,2) -- (42,0);
\draw[orange, very thick] (42,0) .. controls (43,3) and (47,3) .. (48,0) node[black,xshift=-24,yshift=6] {\large $\gamma$};

\draw[orange, very thick] (32.5,-9) .. controls (33.5,-6) and (37.5,-6) .. (38.5,-9) node[black,xshift=-24,yshift=6] {\large $\phi(\beta)$};

\draw[orange, very thick] (39.5,-10) .. controls (40.5,-7) and (44.5,-7) .. (45.5,-10) node[black,xshift=-24,yshift=6] {\large $\phi(\gamma)$};
\draw[very thick] (31.5,-10) -- (32.5,-9);
\draw[very thick] (38.5,-9) -- (39.5,-10);
\draw[gray, thin, dashed] (32.5,-9) -- (38.5,-9);

\draw[very thick,|->] (1,-1.7) -- (1,-5.5);
\draw[very thick,|->] (11,-1.7) -- (11,-5.5);
\draw[very thick,|->] (23.5,-1.7) -- (23.5,-5.5);
\draw[very thick,|->] (38.5,-1.7) -- (38.5,-5.5);
\end{tikzpicture}
\caption{Illustration of the map $\phi$ according to Definition~\ref{Def2}.}
\label{fig4}
\end{figure}

\begin{defn} We  recursively define the map $\phi$ from $\mathcal{S}$ to $\mathcal{D}$ as follows. For $\alpha\in\mathcal{S}$, we set:
$$
\phi(\alpha)=\left\{\begin{array}{llr}
UD&\text{if }\alpha=\epsilon,&(i)\\
U\phi(\beta)D&\text{if }\alpha=E\bar{E}\beta \text{ with } \beta\in\mathcal{S},&(ii)\\
UD\phi(\beta)&\text{if }\alpha=N\bar{N}\beta \text{ with } \beta\in\mathcal{S},&(iii)\\
U\phi(\beta)D\phi(\gamma)&\text{if }\alpha=N\bar{E}\beta E\bar{N}\gamma \text{ with } \beta,\gamma \in\mathcal{S}. &(iv)
\end{array}\right.
$$
\label{Def2}
\end{defn}
Due to the recursive definition, the image of a zigzag knight's path of length $2n$ under $\phi$ is a Dyck path of length $2(n+1)$. For instance, the images of $N\bar{N}$, $E\bar{E}$ and $N\bar{N}N\bar{N}$ are $UDUD$, $UUDD$, and $UDUDUD$, respectively. We refer to Figure~\ref{fig4} for an illustration of this mapping.

The proof of the following theorem is obtain {\it mutatis mutandis} as for Theorem 2 (see Section 2).
\begin{thm}
For all $n\geq 0$, the map $\phi$ induces a bijection between $\mathcal{S}_n$ and $\mathcal{D}_{n+1}$.
\label{bij}
\end{thm}
 For instance, the image of $E\bar{E}N\bar{N}N\bar{E}N\bar{N}E\bar{N}$ under $\phi$ is $UUDUUDUDDUDD$.

 \paragraph{Acknowledgments.}
The authors would like to thank the anonymous referees for useful remarks and comments.

\end{document}